\documentclass[a4paper]{article}
\usepackage{amsmath,amsthm,amsfonts,amssymb}
\usepackage{bm}
\usepackage{hyperref}
\usepackage{algorithm2e}

\newtheorem{theorem}{Theorem}
\newtheorem{dfn}{Definition}
\newtheorem{notation}{Notation}

\providecommand{\Z}{\mathbb{Z}}
\providecommand{\R}{\mathbb{R}}
\providecommand{\p}{\mathcal{P}}

\providecommand{\email}[1]{\href{mailto:#1}{#1}}
\renewcommand{\vec}[1]{\bm{#1}}

\title{A parallel approach to bi-objective integer programming}
\author{W. Pettersson \\
School of Science, RMIT University, \\ 
  Victoria~3000, \textsc{Australia}. \\
  \email{william@ewpettersson.se} \and
M. Ozlen \\
School of Science, RMIT University, \\ 
  Victoria~3000, \textsc{Australia}. \\
  \email{melih.ozlen@rmit.edu.au}}

\date{16 January 2017}

\begin{document}

\maketitle
\begin{abstract}
  To obtain a better understanding of the trade-offs between various objectives, 
  Bi-Objective Integer Programming (BOIP) algorithms calculate the set of all non-dominated vectors and present these as the solution to a BOIP problem.
  Historically, these algorithms have been compared in terms of the number of single-objective IPs solved and total CPU time taken to produce the solution to a problem.
  This is equitable, as researchers can often have access to widely differing amounts of computing power.
  However, the real world has recently seen a large uptake of multi-core processors in computers, laptops, tablets and even mobile phones.
  With this in mind, we look at how to best utilise parallel processing to improve the elapsed time of optimisation algorithms.
  We present two methods of parallelising the recursive algorithm presented by Ozlen, Burton and MacRae.
  Both new methods utilise two threads and improve running times.
  One of the new methods, the Meeting algorithm, halves running time to achieve near-perfect parallelisation.
  The results are compared with the efficiency of parallelisation within the commercial IP solver IBM ILOG CPLEX, and the new methods are both shown to perform better. 
\end{abstract}

\section{Introduction}
Integer Programming (IP) requires either one single measurable objective, or a pre-existing and known mathematical relationship between multiple objectives.
If such a relationship, often called a ``utility function'', is known then one can optimise this utility function (see e.g., \cite{Abbas20061140,Jorge200998}).
However, if the utility function is unknown, we instead identify the complete set of non-dominated solutions for the Bi-Objective Integer Programming (BOIP) problem.
The net result is that a decision maker can more easily see the trade-offs between different options, and therefore make a well informed decision.

Algorithms that determine this complete set can take exact approaches \cite{Benson1998Outer}, or utilise evolutionary techniques \cite{Parsopoulos2002ParticleSwarm}.
For further details on multi objective optimisation see \cite{ehrgott2005multicriteria}.

The performance of BOIP algorithms, and algorithms in general, is often based on the CPU time taken to find the solution.
This does allow algorithms to be compared without needing to use expensive or speciality hardware (as long as comparisons are made on identical hardware setups), but does not take into account real world computing scenarios.
Recent times have seen desktops and laptops move completely to hardware with multiple computing cores.
Given this, we look at how to best utilise parallel processing in BOIP algorithms.
We constrain ourselves to biobjective problems to demonstrate the feasibility of our approach.

In this paper we look at a BOIP algorithm from Ozlen, Burton and MacRae \cite{Ozlen2014moipaira} which calculates a solution by recursively solving constrained IPs.
This algorithm uses IBM ILOG CPLEX to solve all constrained IPs.
As CPLEX does include it's own parallelisation code, we give results on the effectiveness of parallelisation within CPLEX.
We also demonstrate our own approach to parallelisation of the recursive algorithm of Ozlen, Burton and MacRae.
Our new parallel algorithm achieves near-ideal parallelisation, calculating a solution in half the running time without incurring any additional computational time.
This proves far more effective than parallelisation within CPLEX.

Section \ref{sec:background} of this paper gives a background in optimisation.
Section \ref{sec:aira} gives a brief outline of the recursive algorithm we build upon.
Section \ref{sec:parallel} contains details of our parallel computing approach. 
Section \ref{sec:timing} details our software implementation and gives running time comparisons between the original algorithm, the original algorithm with CPLEX parallelisation, and our parallelisations.

\section{Background}\label{sec:background}


In an integer programming (IP) problem, we are given a set of variables and a set of linear inequalities called {\em constraints}.
An assignment of an integer value to each variable is called feasible if it satisfies all constraints, and an assignment which does not satisfy all constraints is called infeasible.
The set of all feasible vectors we will call the {\em feasible set}, and can be defined as follows.
\begin{dfn}
The {\em feasible set} of an IP problem is given by \[
  X = \{ \vec x \in \Z^n | A\vec x = b, x_j \geq 0 \mathrm{\ for\ } j \in \{0,1,\ldots,n-1\} \}
\]
where $A$ is an $n$-by-$n$ matrix and $b$ is an $n$-by-$1$ matrix that together represent the linear constraints of the problem.
\end{dfn}
Note that inequalities can be converted to equalities with the introduction of slack variables (see e.g. \cite{papadimitriou1982combinatorial} or any introductory linear programming text).

Given a feasible set $X$ and an objective function $f$, the goal of an IP is to find the solution $\vec x\in X$ that optimises $f(x)$.
In this paper we assume that all objective functions are to be minimised.
In some scenarios the goal is to maximise a given objective function.
Such a problem can trivially be converted into an equivalent problem where the objective is to be minimised.
We will denote an IP (and various derived problems) with $\p$.


In a multi objective integer programming problem, we do not have one single objective function but rather a set of objective functions.
The goal then is to determine all {\em non-dominated} (or {\em Pareto optimal}) solutions.

\begin{dfn}
Given a pair $f_1,f_2$ of objective functions $f_i: X \rightarrow \R$,
a solution $\vec x \in X$ is considered {\em non-dominated} (or {\em Pareto optimal}) if there does not exist an $\vec{x'} \in X$ with $\vec {x'} \neq \vec {x}$ such that $f_i(\vec{x'}) \leq f_i(\vec{x})$ for all $i \in \{1,2\}$.
\end{dfn}

A {\em bi-objective integer programming} problem then involves the calculation of the set of all non-dominated feasible solutions.



Given a set of objective functions, one of the simpler methods of generating a related single-objective IP is to apply an ordering to the objective functions, and compare solution vectors by considering each objective in order.
That is, each objective function is considered in turn, with objective functions that appear earlier in the ordering being given a high priority.
We will call such a problem a {\em lexicographic bi-objective integer programming} (LBOIP) problem on $k$ objectives.

      

Our parallel algorithm will use different orderings of a set of objective functions to determine the solution set.
To aid readability, we therefore introduce the following notation to refer to different lexicographic variants of a BOIP problem with a given set of objectives.
\begin{notation}
If a lexicographic version of the BOIP $\p$ will order objectives according to the ordered set $(f_1,f_2)$, we will write $\p^{(1,2)}$.
\end{notation}

The optimal solution for a LBOIP will be part of the solution set for the related BOIP, but there is no guarantee that it will be the only solution for the BOIP.
Indeed, it will often not be the only non-dominated solution.
To determine all non-dominated solutions, the algorithm described in Section \ref{sec:aira} utilises {\em constrained lexicographic multi bi-objective linear programming} (CLBOIP).
A CLBOIP is simply a LBOIP with a constraint on the last objective function.
These constraints limit the solution space to some given bound.
\begin{notation}
  Given an LBOIP $\p^{\vec s}$, if the upper bound on the final objective is $l_k$ we will denote the CLBOIP by $\p^{\vec s}( < l_k)$.
\end{notation}
Recall that $\p^{\vec s}$ denotes that the ordering of the objectives is given by the permutation $\vec s$, which will be an element of the symmetric group $S_2$ for our biobjective problems.

%

\section{The algorithm of Ozlen, Burton and MacRae}\label{sec:aira}

The full recursive algorithm as given by Ozlen, Burton and MacRae is suitable for problems with an arbitrary number of objective functions.
Here we give a brief outline of a biobjective version of the algorithm.
For a complete introduction to the algorithm, see \cite{Ozlen2014moipaira}.
Given a biobjective IP problem with two objective functions $f_1$ and $f_2$, the algorithm runs as follows:

\begin{algorithm}[H]\label{algo:orig}
  \KwData{A BOIP $\p$ with objective functions $f_1$ and $f_2$}
  \KwResult{The non-dominated solutions to the BOIP}
  Let $S = \{\}$ be an empty set to which we will add all non-dominated solutions\;
  Let $l_2 = \infty$\;
  \While{$\p^{(1,2)}(<l_2)$ is feasible}{
    Let $\vec x = (x_1, x_2)$ be the optimal vector for the CLBOIP $\p^{(1,2)}(<l_2)$\;
    Add $\vec x$ to $S$\;
    Set $l_2 = x_2$\;
  }
  \caption{A simple overview of the biobjective version of the algorithm from Ozlen, Burton and MacRae.}
\end{algorithm}

The correctness of this algorithm is readily shown by induction.
For a formal proof of the correctness of this algorithm, see \cite{Ozlen2014moipaira}.

\section{Parallelisation}\label{sec:parallel}
Comparisons of algorithms are often based on either the number of single-objective IPs solved, or CPU time taken to solve the problem.
In the real world, the only thing that really matters (for correct algorithms) is the time between describing the scenario and receiving the solution.
Not everyone will have access to supercomputing facilities, however desktop and laptop computers have had multiple cores as standard for many years now.

We therefore look to reduce the elapsed running time of optimisation algorithms by introducing parallelisation.
We will use the term {\em thread} to denote a single computational core performing a sequence of calculations.
In a parallel algorithm, we therefore have multiple {\em threads} which are performing multiple calculations simultaneously.
In this paper we look at improvements gained by utilising two threads at once.

\subsection{Range splitting}
When solving $\p$, it is clear that the maximum and minimum values of $f_1(\vec x)$ can be determined by solving two LBOIP problems $\p^{(1,2)}$ and $\p^{(2,1)}$.
One na\"ive method of distributing this problem across $t$ threads would be to split this range into $t$ equal sized pieces, and then adding an upper and lower bounds on $f_1$ to the specific sub-problem solved by each thread.
These results can be combined in the obvious manner to give the solution.
We will refer to this algorithm as the {\em Splitting} algorithm, as the range of $f_1$ is split up so that each thread gets a single section.
The proof of correctness of this algorithm is trivial.
Implementation and timing results are detailed in Section \ref{sec:timing}.

\subsection{Efficient parallelisation}

Whilst the algorithm discussed in the previous section is parallel, there is no guarantee that all threads will perform an equal (or near-equal) amount of work.
Indeed, it is easy to visualise problems where one thread will perform far more work than another.
Instead we use an algorithm which dynamically adapts itself as the solution set is found.

Recall that in Algorithm \ref{algo:orig} we used the specific ordering $\p^{(1,2)}$.
We could also solve $\p^{(2,1)}$ and obtain the same result.
This idea forms the basis of our work.
We show below how the limit $l_2$ obtained from $\p^{(1,2)}(< l_2)$ is able to be shared with the problem $\p^{(2,1)}(< l_1)$.
This allows the two problems to be solved simultaneously, which almost halves the running time of our new algorithm when compared to the original.

\begin{theorem}\label{thm:parallel}
  If we have the complete set $S$ of non-dominated solutions for $\p^{(1,2)}$ with $x_2 \geq k$,
  the complete set $S'$ of non-dominated solutions for $\p^{(2,1)}$ with $x_1 \geq l$, and we also have the non-dominated solution $(l,k)$, then the union $S \bigcup S' \bigcup \;\{(l,k)\}$ is the complete set of non-dominated solutions to $\p$.
\end{theorem}

\begin{proof}
Assume that $\vec x = (x_1,x_2)$ is a non-dominated solution to $\p$ that is not in either $S$ nor $S'$.
If $x_1 > l$ then $\vec x \in S'$, a contradiction.
Similarly, if $x_2 > k$ then $\vec x \in S$ which is also a contradiction.
Therefore $x_1 \leq l$ and $x_2 \leq k$, but then as $(l,k)$ is non-dominated, the only solution is $(x_1,x_2) = (l,k)$.
\end{proof}

Note that $\p^{(1,2)}$ and $\p^{(2,1)}$ are both CLBOIP problems that can be solved independently by Algorithm \ref{algo:orig}, and that $(l,k)$ will be found as a solution to both of these problems.
Given this result, we propose the following parallel algorithm for computing the solution to BOIP problems.

%
\begin{algorithm}[H]\label{algo:new}
  \KwData{A BOIP $\p$ with objective functions $f_1$ and $f_2$}
  \KwResult{The non-dominated solutions to the BOIP}
  Let $t \in \{1,2\}$, and let $t'$ be the unique value in $\{1,2\}\setminus\{t\}$\;
  Let $s_1 = (2,1)$ and $s_2 = (1,2)$\;
  Let $S_1 = S_2 = \{\}$ be empty sets\;
  Let $l_1 = l_2 = \infty$\;
  \ForEach{thread $t$}{
    \While{$\p^{s_t}(<l_t)$ is feasible}{
      Let $\vec x = (x_1, x_2)$ be the solution for the CLBOIP $\p^{s_t}(<l_t)$\;
      Add $\vec x$ to $S_t$\;
      Set $l_t = x_t$ *\;
      Add $x_{t'} < l_{t'}$ as a constraint to $\p^{s_t}(< l_t)$ *\;
    }
  }
  \Return {$S_1 \cup S_2$}
  \caption{A parallel version of the algorithm from \cite{Ozlen2014moipaira}.
  Note that in the lines marked *, the values $l_{t}$ and $l_{t'}$ are shared between the two threads.}
\end{algorithm}
We call this algorithm the {\em Meeting} algorithm, as the two threads meet in the middle to complete the calculations.
Correctness of the Meeting algorithm follows from Theorem \ref{thm:parallel} and the correctness of Algorithm \ref{algo:orig}.

\section{Implementations and running times}\label{sec:timing}

\subsection{Implementation}
Our algorithms were implemented in C++, and are available at \\
\url{https://github.com/WPettersson/moip\_aira}.
All calculations were run on the NCI supercomputing cluster Raijin, on nodes consisting of two Intel Sandy Bridge E5 2670 processors and 32GB of RAM. Code was compiled with GCC 4.9, using no special optimisation controls beyond \texttt{-O2}.
We compared the elapsed running time (and not computational time) of the original algorithm (with both one thread allocated to CPLEX, and two threads allocated to CPLEX), along with the {\em Splitting} algorithm and the {\em Meeting} algorithm.
These running times are summarised in Table 1.
\begin{table}
  \begin{center}
\begin{tabular}{|c|c|c|c|c|}
  \hline
  \multicolumn{5}{|c|}{Assignment problems} \\
  \hline
  \# tasks & Ozlen et al. & CPLEX & Splitting & Meeting \\
  \hline
  40 & 10.95 & 11.00 & 9.14 & 5.74 \\
  60 & 34.42 & 31.89 & 28.74 & 17.83 \\
  80 & 68.39 & 57.55 & 55.57 & 35.63 \\
  100 & 118.69 & 106.30 & 95.66 & 63.37 \\
  200 & 515.57 & 453.54 & 402.98 & 276.90 \\
  500 & 3262.26 & 3468.03 & 2327.63 & 1738.74 \\
  \hline
  \multicolumn{5}{|c|}{Knapsack problems} \\
  \hline
  \# items & Ozlen et al. & CPLEX & Splitting & Meeting \\
  \hline
  50 & 1.00 & 1.10 & 0.67 & 0.53 \\
  100 & 5.03 & 4.83 & 3.60 & 2.59 \\
  200 & 22.37 & 20.56 & 16.13 & 11.53 \\
  400 & 73.75 & 71.69 & 57.70 & 36.42 \\
  1000 & 338.67 & 347.34 & 263.01 & 150.06 \\
  2000 & 1200.50 & 1113.11 & 912.85 & 528.85 \\
  \hline
\end{tabular}
\end{center}
\label{tab:problems}
\caption{Elapsed running timing comparisons of the four algorithms. We ran ten different random versions of each sized problem and averaged the results.}
\end{table}

From the table, we see that letting CPLEX use a second thread improved running times only slightly.
It is not surprising that CPLEX does not parallelise efficiently in this manner, as CPLEX cannot take advantage of the full details of the algorithm used.
The Splitting algorithm was more impressive, showing significant results.

However, our Meeting algorithm is the clear outlier, twice as fast as the original algorithm of Ozlen, Burton and MacRae on all problems.
This is the ideal outcome for parallelising the algorithm with two threads.

\section{Conclusion}
We successfully implemented two parallel algorithms to solve biobjective optimisation problems.
Both improved performance, with one showin ideal performance increase.
For biobjective problems (and potentially multi objective problems) this faster algorithm allows solutions to be found in half the time.
Ongoing work in this field will look at various ways of utilising more threads in parallel to further improve running times for IP problems with three or more objectives.

\paragraph{Acknowledgements}
Melih Ozlen is supported by the Australian Research Council under the Discovery Projects funding scheme (project DP140104246).

\bibliography{bibliography}
\bibliographystyle{plain}
\end{document}